\documentclass[reqno,a4paper, 11pt]{amsart}

\usepackage{enumerate}
\usepackage{amsmath,amssymb,amscd,amsfonts,verbatim}
\usepackage[arrow,matrix]{xy}
\usepackage{amsmath}

\newtheorem{thm}{Theorem}[section]
\newtheorem{lem}[thm]{Lemma}
\newtheorem{prop}[thm]{Proposition}
\newtheorem{cor}[thm]{Corollary}

\theoremstyle{definition}
\newtheorem{rem}[thm]{Remark}

\newtheorem{assu}[thm]{Assumption}

\newcommand{\Pic}{{\text{\rm Pic}}}
\newcommand{\Tors}{{\text{\rm Tors}}}

\newcommand{\Aut}{\text{\rm Aut}}

\newcommand{\OO}{{\mathcal {O}}}

\newcommand{\inv}{^{-1}}
\newcommand{\al}{\alpha}

\newcommand{\Ga}{\Gamma}

\newcommand{\De}{\Delta}
\newcommand{\epsi}{\epsilon}

\newcommand{\fie}{\varphi}
\newcommand{\si}{\sigma}
\newcommand{\Si}{\Sigma}

\newcommand{\C}{{\mathbb C}}

\newcommand{\Q}{{\mathbb Q}}
\newcommand{\Z}{{\mathbb Z}}
\newcommand{\F}{{\mathbb F}}
\newcommand{\pp}{{\mathbb{P}}}
\newcommand{\tpi}{\tilde{\pi}}
\newcommand{\tsi}{\tilde{\si}}

\numberwithin{equation}{section}
\title[Involutions on numerical Campedelli surfaces]
{Involutions on numerical Campedelli surfaces}
\author{Alberto Calabri, Margarida Mendes Lopes and Rita Pardini}

\email{calabri@dm.unibo.it}
\curraddr{Dipartimento di Matematica, Universit\`a degli Studi di Bologna\\
Piazza di Porta San Donato, 5 --- I-40126 Bologna \\Italy}

\email{mmlopes@math.ist.utl.pt} \curraddr{Departamento de
Matem\'{a}tica, Instituto Superior T\'ecnico,
Universidade T{\'e}cnica de Lisboa\\
Av.~Rovisco Pais ---  P-1049-001 Lisboa \\ Portugal}

\email{pardini@dm.unipi.it} 
\curraddr{Dipartimento di
Matematica, Universit\`a degli Studi di Pisa \\ Largo B. Pontecorvo 5,
56127  Pisa \\Italy}

\thanks{{\it Mathematics Subject Classification (2000)}: 14J29. }

\begin{document}
\begin{abstract}

 Numerical Campedelli surfaces are minimal surfaces of general type
 with $p_g=0$ (and so $q=0$) and $K^2=2$. Although they have been studied by several authors,
their complete classification is not known.

In this paper we classify numerical Campedelli  surfaces
with an involution, i.e.\ an automorphism of order 2. First  we show that an 
involution on a numerical  Campedelli surface $S$ has either   four or six 
isolated 
fixed 
points, and the  bicanonical map of $S$ is composed with the involution
 if and  only if the involution has six isolated fixed points. Then we study in detail each of the possible cases, describing also several examples.

\end{abstract}
\maketitle

\section{ Introduction}

 Numerical Campedelli surfaces are minimal surfaces of general type
 with $p_g=0$ (and so $q=0$) and $K^2=2$. The first such example was 
presented by 
Campedelli (\cite{Cam}) in 1932. Since then several authors
 (cf. \cite{miyaoka1}, \cite{peters}, \cite{reid}, \cite{Re2}, 
\cite{kotschick}, \cite{Su}, \cite{Naie2}, \cite{kulikov},...) 
 have studied these 
surfaces, but our knowledge about them is far from being complete.

Since a classification  of numerical Campedelli surfaces does not
 seem feasible at the moment, a possible approach is
 to restrict one's attention to the Campedelli surfaces
which have  some additional geometrical feature. This is
 what we do in the present paper, where we study the Campedelli surfaces
 which have an involution, i.e. which have  an automorphism of order 2. This choice  is motivated by work of Keum and Lee (\cite{keum}) and  of Calabri, Ciliberto and Mendes Lopes (\cite{CCM2}), who have studied the same problem for numerical Godeaux surfaces, 
that is  minimal surfaces of general type with $p_g=0$
 and $K^2=1$.
 
  In order to put our work  in perspective, we briefly recall here the main results of the paper \cite{CCM2}, which contains  a complete classification 
  of numerical Godeaux surfaces  with an involution.
  
  If $S$ is  a numerical Godeaux surface and $\si$ is an involution  of $S$, then $\si$ has five isolated fixed points and:
  \begin{itemize}
  \item the 
bicanonical map of the surface factors through the natural projection onto 
the quotient surface $S/\si$;
\item the quotient surface is either rational or birational to an Enriques surface;
\item the possible quotient surfaces are classified and examples of  each possibility in the list do exist;
 \item if $S/\si$ is rational,   then the surface $S$ can be   obtained  as a specialization of one of the surfaces in  the list proposed by 
Duval (see \cite{SantaCruz}, cf. also \cite{borrelli}), by letting the 
branch locus acquire   some 
singularities.
  \end{itemize}
  
In the case of numerical Campedelli surfaces the situation is more involved, since the bicanonical map  may not factor through the quotient map $S\to S/\si$. 
Indeed, we show that an 
involution on a numerical  Campedelli surface $S$ has either   four or six 
isolated 
fixed 
points, and the  bicanonical map of $S$ factors through the quotient map $S\to S/\si$ if and  only if the involution has six isolated fixed points. In the latter case the
situation is very similar to the case of Godeaux surfaces.  We have the following:
\begin{itemize}
\item the
ramification divisor $R$ on $S$ is not 0, and its components
  can be described (see \S 3);
\item   the quotient surface $S/\si$ 
is either birational to an Enriques surface or a rational surface;
\item if $S/\si$ is rational, then there are four possible cases which 
all 
have a precise description (see also \S 3). Each of the four cases actually occurs (cf. \S 5). 
\end{itemize}
The analysis in \S 3 shows also that, if the bicanonical map of $S$ is composed with the involution, then   the $2-$torsion of the  surface $S$ is non trivial in three of the possible five cases.

If the bicanonical map is not composed with the involution, i.e. if the 
involution has four isolated fixed points,  we show that the ramification 
divisor $R$ 
is either $0$ or constituted by one, two or three $-2-$curves. Note that 
if $R\neq 0$ then $K_S$ is not ample.

In this case there are more possibilities for the quotient surface $S/\si$, as explained below:
\begin{itemize}
\item $S/\si$ is of general type (a numerical Godeaux surface) if and only 
if the ramification divisor $R$ is equal to 0;
\item if $R$ is irreducible, then $S/\si$ is properly elliptic;
\item if $R$ has two or three components then $S/\si$ may be rational or
birational to an Enriques surface or properly elliptic. 
\end{itemize}
The  case where $S/\si$ is a numerical  Godeaux surface appears  in the 
examples constructed by R. Barlow in \cite{Ba1} and \cite{Ba2}.
 In \S  5, by specializing one of the examples of R. 
Barlow, we present examples for which the quotient surface is either an 
elliptic surface or birational to 
an Enriques surface. We do not know any instance in which the quotient 
surface  is a 
rational surface for this case.

 In \S 5 we also study a family 
of numerical 
Campedelli surfaces with torsion $\Z_3^2$, whose construction is 
attributed by J. Keum to   A. 
Beauville and X. Gang. We show that every surface in this family  
has two involutions, one 
with four isolated fixed points and one with six isolated fixed points, 
whose quotients are respectively birational to a numerical Godeaux surface 
and a 
rational surface.

In \S 5 we study the involutions of numerical Campedelli surfaces
 with torsion $\Z_2^3$,  the so-called ``classical Campedelli surfaces''. 
Using the description of these surfaces as a $\Z_2^3$-cover of $\pp^2$ 
branched on 7 lines (cf. \cite{kulikov}), 
 we show that these involutions are all composed with the bicanonical map.

The paper is organized as follows.  In \S 2, using the results in 
\cite{CCM2}, we describe the general properties of numerical Campedelli 
surfaces with an involution, showing in particular that such an involution 
always has four or six isolated fixed points.

In \S 3 we study the case where the involution has six isolated fixed 
points and we describe in some detail each possibility.
In \S 4 we study the case where the involution has four fixed 
points. Finally in \S 5 we describe the examples, two of which were not (to our knowledge) previously known.

\medskip 

{\bf Acknowledgments.} The original idea for this work is due to discussions with Ciro Ciliberto, who should be also considered  an author of this work, and whom we wish to thank warmly. We are also indebted to JongHae Keum for  communicating to us Example 2 of \S 5. 

The first and the last author are members of G.N.S.A.G.A.-I.N.d.A.M.
 The second author is a member of the Center for Mathematical
Analysis, Geometry and Dynamical Systems. This research was partially supported by the italian  project "Geometria sulle
variet\'a algebriche" (PRIN COFIN 2004) and by FCT (Portugal) through program POCTI/FEDER and Project
POCTI/MAT/44068/2002. 

The first author is very grateful to Mirella Manaresi, and
to all the Algebraic Geometry Group in Bologna,
for providing a stimulating atmosphere
and supporting his research.

\medskip

{\bf Notation and conventions.} We work over the complex numbers. All varieties are projective.

Most of the notation is standard in algebraic geometry, hence we only recall here  a few conventions that we use and that  are   maybe  not universally accepted. 
We denote linear equivalence of divisors on a smooth variety by $\equiv$ and numerical equivalence by $\sim$.  A divisor $D$ on a smooth variety  $X$ is said to be {\em even} if its class is divisible by 2 in the group $\Pic(X)$. 

An {\em involution} of a variety is a biregular automorphism of order 2.  A map $f\colon X\to Y$ of projective varieties  is said to be {\em composed with} an involution $\si$ if $f\circ \si=f$. 
A {\em $-n-$curve}  on a smooth surface is a curve $C$ such that $C\simeq\pp^1$ and $C^2=-n$.

\section{Involutions on a numerical Campedelli surface}\label{prima}

Throughout all the paper we make the  following:
\begin{assu} $S$ is a smooth minimal complex projective surface of general type with  $p_g(S)=0$, $K^2_S=2$ (hence also $q(S)=0$).
Such a surface $S$ is called a {\em numerical Campedelli surface}. 

Moreover we assume that we are given an involution $\si$ of $S$, namely an automorphism $\si\colon S\to S$ of order 2. 
\end{assu}

In this section we establish the notation and recall  some known facts on involutions,  giving all the statements  in the special case of a numerical Campedelli surface. Our main reference is the paper \cite{CCM2}, which contains a detailed analysis of involutions on surfaces of general type with $p_g=0$.
\bigskip

The fixed locus of the involution $\si$ is the union of
an effective divisor $R$ and of $k$ isolated points $p_1,\ldots,p_k$. The effective divisor $R$, if not $0$, is a smooth, possibly reducible, curve. 
Let $\pi\colon S\to \Si:=S/\si$ be the quotient map,
and set $B:=\pi(R)$ and $q_i:=\pi(p_i)$,
$i=1,\ldots,k$. The surface $\Si$ is normal and $q_1,\ldots,q_k$  are
ordinary double points, which are the only singularities of $\Si$.
In particular, the singularities of $\Si$ are canonical and the
adjunction formula gives $K_S\equiv\pi^*K_{\Si}+R$.

Let $\epsi\colon V\to S$ be the blow-up of $S$ at $p_1,\ldots,p_k$
and let $E_i$ be the exceptional curve over $p_i$, $i=1,\ldots,k$.
Then $\si$ induces an involution $\tsi$ of $V$ whose fixed locus
is the union of $R_0:=\epsi^*(R)$ and of $E_1,\ldots,E_k$.
Denote by
$\tpi\colon V\to W:=V/\tsi$ the projection onto the quotient and
set $B_0:=\tpi(R_0)$, $N_i:=\tpi(E_i)$, $i=1,\ldots,k$. The surface
$W$ is smooth and the $N_i$ are disjoint $-2-$curves.
Denote by $\eta\colon W\to \Si$ the map induced by $\epsi$.
The map $\eta$  is the minimal resolution of
the singularities of $\Si$ and there is a commutative diagram:
\begin{equation}\label{diagram}
\begin{CD} V@>{\epsi}>>S\\ @V{\tpi}VV @VV{\pi}V\\
W@>{\eta}>>\Si
\end{CD}
\end{equation}
The map $\tpi$ is a flat double cover branched on
$\tilde B:=B_0+\sum_{i=1}^k N_i$,
hence there exists a divisor $L$ on $W$
such that $2L\equiv \tilde B$, namely $\tilde{B}$ is an even divisor.
\begin{rem}\label{Wregular}
We have $p_g(V)=q(V)=0$, since $V$ is birational  to $S$. Since $V$ dominates $W$, we also have $p_g(W)=q(W)=0$.
\end{rem}

The number $k$ of isolated fixed points  is a very important invariant of the involution $\si$. As explained below, it determines whether the bicanonical map $\fie\colon S\to\pp^2$ is composed with $\si$.
\begin{prop}[\cite{CCM2}, Proposition 3.3, (v) and Corollary 3.6]\label{cases}
One of the following two possibilities occurs:
\begin{itemize} 
\item[I)] $k=6$. In this case $\fie$ is composed with $\si$.
\item[II)] $k=4$. In this case $\fie$ is not composed with $\si$. More precisely,   $\pi^*|2K_{\Si}+B|$ has dimension 1, namely it is a codimension 1 subsystem of $|2K_S|$.
\end{itemize}
\end{prop}

We set   $D:=2K_W+B_0$. The divisor $D$  will play an important role in our analysis of numerical Campedelli surfaces with an involution.

One has the following properties (cf. \cite[\S 3]{CCM2},  for the proofs):
\begin{prop}\label{properties}
\begin{enumerate}[\rm(i)]
%\item $K_V=\tpi^*(K_W+L)$; 
\item $\epsi^*(2K_S)=\tilde{\pi}^*D$;
\item $D$ is nef and big,  and $D^2=4$;
%\item $H^i(W,\OO_W(2K_W+L))=0$, $i>0$.
\item $D+N_1+\dots +N_k$ is an even divisor;
%\item $K_W L+L^2=-2$;
%\item $K^2_W+K_W L\ge 0$;
%\item $h^0(W, \OO_W(2K_W+L))=K_W^2+K_W L$;
\item if $k=6$, then: 
$-4 \le K^2_W\le 0,\ \  K_WD=0$;
 \item if $k=4$, then:
 $-2\le K^2_W\le 1, \ \ K_WD=2$.
 \end{enumerate}
\end{prop}
\begin{rem}\label{KC}
We will often apply Proposition \ref{properties}, (i) as follows. Given a curve $C$ of $W$, we can pull it back to a curve $C'$ of $V$. If $C'$ is not contained in the exceptional locus of $\epsi$, then we  can push it down to a curve $\tilde{C}$ on $S$. Then $K_S\tilde{C}=DC$.
\end{rem}

%%%%%%%%%%%%%%%%%%%%%%%%%%%%
Assume  that $K_W+D$ is not nef.  Then one can show that there is an irreducible $-1-$curve $E$ on $W$ with $DE=0$, $E N_i=0$ for $i=1,\dots, k$. By repeatedly blowing down such $-1-$curves, one obtains a sort of minimal model for the pair $(W, K_W+D)$.
More precisely, we have the following:
\begin{prop}[\cite{CCM2}, Proposition 3.9]\label{W'}
There exists a birational morphism $f\colon W\to
W'$, where $W'$ is smooth,  with the following properties:
\begin{enumerate} [\rm (i)]
\item  for $i=1,\dots, k$ the curve  $N'_i:=f(N_i)$ is a   $-2-$curve on $W'$ and the curves $N'_1,\dots, N'_k$ are disjoint;
\item there is a nef  divisor $D'$ on $W'$ such that $f^*(D')=D$, ${D'}^2=D^2$ and $K_{W'}D'=K_WD$;
\item $D'N'_i=0$ for $i=1,\dots, k$ and $D'+N'_1+\dots+N'_k$ is an even divisor;
\item
$K_{W'}+D'$ is nef.
\end{enumerate}
\end{prop}
\begin{rem}\label{remW'}
 The proof of \cite{CCM2}, Proposition 3.9 actually shows more. Namely:
\begin{enumerate}[\rm(i)]
\item since $K_W+D$ is effective and the curves contracted by $f$ satisfy $E(K_W+D)<0$, the components of the exceptional locus of $f$ are contained in the fixed part of $|K_W+D|$;
\item if $E$ is an irreducible component of the exceptional locus of $f$, then $E$ gives a $-2-$curve on $S$. In particular, if $K_S$ is ample then we have $W=W'$.
\end{enumerate}

\end{rem}

\section{Involutions composed with the bicanonical map}\label{sCampedelli}
This section is devoted to the study of case I) of Proposition \ref{cases}, namely here we assume that $k=6$ and  the bicanonical map $\fie\colon S\to\pp^2$ is composed with $\si$. 

In what follows we use  freely the notation introduced in section 2.
By Proposition \ref{properties}, in this case $-4\leq K_W^2\leq 0$ and $K_WD=0$. 
This allows us to  establish some properties of the ramification divisor $R$ on $S$.

 Using Proposition \ref{properties}, and arguing as in the proof of Proposition 4.5 of \cite{CCM2}, one obtains the following:

\begin{prop}\label{godeaux}
Let $S$ be a numerical Campedelli surface with an involution $\si$, such that 
  the bicanonical map $\fie$ is composed with $\si$. Then the divisorial part $R$ of the fixed locus of $\si$ satisfies: 
\begin{enumerate} [\rm(i)]

\item $K_S R=2$;
\item $R^2= 2K_W^2+2$ is even,  and $-6 \leq R^2 \leq 2$.

\end{enumerate}

\noindent Furthermore $R=\Gamma+ Z_1+\cdots+Z_h$  
where $\Gamma$ is a smooth curve
with $K_S \Gamma=2$ and $Z_1,\ldots,Z_h$ are disjoint $-2-$curves,
which are disjoint also from $\Gamma$. Here

\begin{enumerate} [\rm(i)] \addtocounter{enumi}{2}

\item either $\Gamma$ is irreducible, $0\le p_a(\Gamma)\le 3$ and
$\Gamma^2=2p_a(\Gamma)-4$;
or $\Gamma$ has exactly two components $\Gamma_1+\Gamma_2$, where $\Gamma_i$, $i=1,2$, is either a 
rational curve
with self-intersection $-3$ or an  elliptic curve  with
self-intersection $-1$;
\item the number $h$ of   $-2-$curves $Z_1,\ldots,Z_h$
satisfies
\[
h=p_a(\Gamma)-K_W^2-3 \ge 0;
\]

\item if $\Gamma^2=2$, then $\Gamma\sim K_S$ and $S$ has non-trivial
torsion.
\end{enumerate}

\end{prop}

\bigskip

In order to study in more detail these surfaces  
 we  consider  the system $|D|:=|2K_W+B_0|$ and  its adjoint 
systems.

\begin{lem}\label{emme} Let $|K_W+D|=|M|+F$, where $F$ is the fixed part.
Then one has:
\begin{enumerate}[\rm(i)]
\item  $h^0(W,\OO_W(M))=3$;

\item  $M D=4$;
\item  if $F\neq 0$, then every component $E$ of $F$ is
such that $D E=0$ and $E^2<0$.
\end{enumerate}
\end{lem}

\begin{proof}
Assertion (i) follows by the adjunction sequence for the general $D$, since $W$ is regular by Remark \ref{Wregular} and $p_a(D)=3$ by Proposition \ref{properties}.

Let us prove part (ii).
Since, by Proposition  \ref{properties},  $D$ is nef and $D (M+F)=4$, one has $M D\leq 4$.

Suppose by contradiction that $MD< 4$. We claim that in this case $|M|$ is not composed with a pencil, and so, in particular, $M^2>0$. Indeed,  if $|M|=|2C|$ and $MD<4$ then one would have  $CD\leq 1$.
But then  $|C|$ would give  a pencil $|\tilde C|$  on $S$ such that $\tilde C K_S=1$ (cf. Remark \ref{KC}), which is impossible by the index theorem. 

By a similar argument we verify that  $MD\geq 2$.

Suppose that $MD=2$.
Then by the index theorem  we obtain $M^2=1$ and $2M\sim D$. This is impossible, because we have  $K_WD=0$ by Proposition \ref{properties} and this implies $K_WM=0$, contradicting  the adjunction formula.

So we are left with the case $MD=3$.
We have $M^2+MF=M(K_W+D)$.  So $MD=3$ means that $MF$ is also odd and thus, because $M$ is nef, $MF\geq 1$. Therefore $MK_W=M^2+MF-3\geq M^2-2$.

 On the other hand, the index theorem gives $M^2D^2\leq (MD)^2= 9$ hence $M^2\leq 2$.
 Since  $|M|$ is not composed with a pencil, we have  $M^2=1 $ or $M^2= 2$.

In the first case $\phi_M\colon W\to \pp^2$ is a birational  morphism,
but this is impossible because  $2g(M)-2=M^2+K_WM\geq 0$.

In the second case the system $|M|$ gives a system $|\tilde{M}|$ on $S$ with $\tilde{M}^2\ge 4$ and $K_S\tilde{M}=3$ (cf. Remark \ref{KC}). By the adjunction formula we get $\tilde{M}^2\ge 5$, contradicting the index theorem applied to $K_S$ and $\tilde{M}$.

So we have shown that $MD=4$. Now (iii) follows immediately from $DF=0$ and the index theorem.
\end{proof}

Consider now the morphism $f\colon W\to W'$ and the divisor $D'$ of  Proposition \ref{W'}.

\begin{prop}\label{K2W'}\begin{enumerate}[\rm(i)]
\item One has:  $-4\le K_{W'}^2 \le 0$;
\item if $K^2_{W'}=0$, then $W'$ is an Enriques surface;
\item  
if $K_{W'}^2<0$, then $W'$ is rational.
\end{enumerate}
\end{prop}

\begin{proof} We recall that $-4\le K^2_W\le 0$ by Proposition \ref{properties}, and so $K^2_{W'}\ge -4$.
Since $D' (K_{W'}+D')=4$, the index theorem implies that
$(K_{W'}+D')^2\leq 4$, or equivalently $K_{W'}^2\le 0$.

The surface $W'$ is either rational or birational to an Enriques surface by \cite[Corollary 3.7]{CCM2}. Since $K_{W'}D=0$,  if $K^2_{W'}=0$ then $K_{W'}\sim 0$ and $W'$ is an Enriques surface. 
If $K_{W'}^2<0$, then $K_{W'}(K_{W'}+D)<0$. Since $K_{W'}+D'$ is nef by Proposition \ref{W'}, this implies that the Kodaira dimension of $W'$ is negative, and therefore $W'$ is rational.
\end{proof}

\begin{lem}\label{fixed}
If $K_{W'}^2<0$, then $|K_{W'}+D'|$ has
no fixed part.
\end{lem}

\begin{proof}
Write, as usual, $|K_{W'}+D'|=F'+|M'|$, where $F'$ is the fixed
part.   Since 
the morphism $f\colon W\to W'$ contracts only
curves that are fixed for $|K_W+D|$ (cf. Remark \ref{remW'}), by Lemma \ref{emme} we see that $M' D'=4$, $F'D'=0$.

Notice that $F' K_{W'}=F'(K_{W'}+D')=F' M'+{F'}^2$, so $F'
M'$ is even. Since both $M'$ and $K_{W'}+D'$ are nef (cf. Proposition \ref{W'}, (iii)), we have the following inequalities:
$$M' F'\leq M' F'+{M'}^2=M'(K_{W'}+D')\leq
(K_{W'}+D')^2=K_{W'}^2+4<4.$$
It follows that $M' F'=0$ or $M'F'=2$.  If $M'F'=2$, then ${M'}^2\leq 1$.

We start by seeing that $M'F'=2$ does not occur.
If $M'F'=2$ and ${M'}^2=1$, then $K_{W'}M'=(K_{W'}+D')M'-4=3-4=-1$. Since ${M'}^2=1$,  $|M'|$ is not composed with a pencil, the general curve of $|M'|$ is smooth and  $\phi_{M'}\colon W'\to \pp^2$ is a birational  morphism. This is impossible because $p_a(M')=1$.

If $M'F'=2$ and ${M'}^2=0$, then $M'\equiv 2C$, where $|C|$ is a free pencil.  Since  $K_{W'}M'=(K_{W'}+D')M'-4=2-4=-2$, one has that $K_{W'}C=-1$, which contradicts the adjunction formula. So $M'F'=2$ does not occur.

 On the other hand, if  $M'F'=0$ then  $F'=0$. Otherwise, since $D' F'=0$, then
${F'}^2<0$, implying that $F' (K_{W'}+D')= {F'}^2+ M'
F'<0$. This contradicts the fact that $K_{W'}+D'$ is nef.
So  $|K_{W'}+D'|$ has no fixed part.
\end{proof}

 Next we examine separately each of the possibilities for $K^2_{W'}$, which ranges between $-4$ and $0$ by Proposition \ref{K2W'}. 

\medskip
\subsection{The case $K^2_{W'}=0$}

In this case the surface $W'$ is an Enriques surface by Proposition \ref{K2W'}. 
\begin{prop}\label{Dfree} The system $|D'|$ is base point free and irreducible.
\end{prop}
\begin{proof} Write $D':=|M|+F$, where $F$ is the fixed part. By Proposition \ref{properties} (i) and Proposition \ref{W'} (ii), the system $|M|$ pulls back on $S$ to the moving part of $|2K_S|$. Since the bicanonical image of $S$ is a surface by \cite{xiaocan}, the general $M$ is irreducible.
In particular, $M$ is nef and big and the Riemann--Roch theorem gives $3=h^0(M)=M^2/2+1$, namely $M^2=4$.
So we have: $4=M^2\le M^2+MF\le D^2=4$, which implies $MF=F^2=0$. Hence $F=0$ by the index theorem.

Now assume that $|D'|$ has base points. By Proposition 4.5.1 of \cite{CD}, there exists an effective divisor $E$ on $W$ such that $E^2=0$, $ED'=1$. By Remark \ref{KC}, this gives a divisor $\tilde{E}$ on $S$ with $K_S\tilde{E}=1$ and $\tilde{E}^2\ge 0$. The adjunction formula then gives $\tilde{E}^2\ge 1$, but this contradicts the index theorem.
\end{proof}

\begin{cor}
The bicanonical system $|2K_S|$ is base point free.
\end{cor}
\begin{proof}
The statement  follows immediately by Proposition \ref{Dfree}, since $|2K_S|$ is the pull back of $|D|$ to $S$ by Proposition \ref{properties}.
\end{proof}

\begin{prop}\label{torsionenriques} The torsion group $\Tors(S)$ of $S$ has order 4 or 8.
\end{prop}
\begin{proof} 
Since the group $\Tors(S)=\Tors(V)$ has order at most  9 (cf. \cite[Chap. VII.10]{bpv2}), it is enough to show the existence of  an \'etale cover of $V$ of degree 4.

Let $p\colon K\to W$ be the \'etale double cover of $W$ induced by the  K3 cover of $W'$. Then we have a cartesian diagram:
\begin{equation*}
\begin{CD}\tilde{K} @>{\tilde{p}}>>V\\
@V\rho VV @VV\tilde{\pi}V\\
K @>{p}>> W
\end{CD}
\end{equation*} 
The map $\tilde{p}$ is an \'etale double cover, while $\rho$ is a double cover branched on the inverse image  $\De$ of $B_0+N_1+\dots +N_6$. The divisor $\De$ is the disjoint union of a  divisor $\De_0$ with $\De_0^2=8$ and of twelve $-2-$curves $\Ga_1,\dots, \Ga_{12}$.
Consider the natural map $\psi\colon \Z\De_0\oplus\Z\Ga_1\oplus\dots \oplus \Z\Ga_{12}\to H^2(K, \Z_2)$. The image of $\psi$ is a totally isotropic subspace, hence it has dimension at most 11, since $h^2(K,\Z_2)=22$ and the intersection form on $H^2(K, \Z_2)$ is non degenerate by Poincar\'e duality. Hence the kernel of $\psi$ has dimension at least 2. By Lemme 2 of \cite{beauville}, the surface $\tilde{K}$ has a connected \'etale double cover, hence $V$ has a connected \'etale cover of degree 4.
\end{proof}

\begin{rem}  Examples of this situation can be found in \cite{Naie}.
Those examples have torsion group $\Z_2^3$ or $\Z_2\times \Z_4$.
\end{rem}
\medskip 

\subsection{The case $K_{W'}^2=-1$} \label{K2W'=-1}

By Proposition \ref{K2W'}, $W'$ is a rational surface. Denote  $M':=K_{W'}+D'$ and recall that  $|M'|$ has no fixed part by Lemma \ref{fixed}. One has ${M'}^2=3$, 
  $K_{W'}M'=-1$. Since $|M'|$ is 2-dimensional,  the general curve of $|M'|$ is irreducible. The system $|K_{W'}+M'|$ has dimension 1 by the adjunction sequence for the general $M'$. 
\begin{lem} \label{C} The linear system $|K_{W'}+M'|$ is a base point free pencil of non-singular rational curves.
\end{lem}
\begin{proof}
We claim that $K_{W'}+M'$ is nef. 

Suppose otherwise. Then there exists an irreducible curve $\theta$ with $\theta(K_{W'}+M')<0$. It follows that $\theta$ is a fixed component of $|K_{W'}+M'|$ and $\theta^2<0$. The general $M'$ is smooth and irreducible and $|K_{W'}+M'|$ restricts to the complete canonical system on $M'$. Hence the general $M'$ does not meet $\theta$, namely $\theta M'=0$ and $\theta K_{W'}<0$. Thus $\theta$ is a $-1-$curve.

The divisor  $G:=M'-\theta$ is effective, since $|M'|$ has dimension 2, and we have $G^2=2$, $GD'=3$. Then $G$ gives a divisor $\tilde G$ on $S$ such that $\tilde G^2\geq 4$, $\tilde  G K_S=3$ (cf. Remark \ref{KC}), and therefore $\tilde G^2\ge 5$ by the adjunction formula. This contradicts the index theorem applied to $K_S$ and $\tilde G$, showing that $K_{W'}+M'$ is nef.

Consider $|K_{W'}+M'|= |C|+F $, where $F$ is the fixed part. Arguing as above, one shows that 
$M'F=0$,  and therefore $F^2<0$ if $F\neq 0$. Because $K_{W'}+M'$ is nef, $C(K_{W'}+M')=C^2+CF\geq 0$ and $F(K_{W'}+M')=F^2+CF\geq 0$. Since $(K_{W'}+M')^2=0$,  we have equality in both cases. 

But then, because $C$ is nef, we must have $CF=0$,  implying also $F^2=0$ and so $F= 0$.
So $|K_{W'}+M'|=|C|$ is a pencil of rational curves.
\end{proof}

\begin{prop}\begin{enumerate}[\rm(i)]
\item There exists a fibration $f\colon S\to\pp^1$ with 3 double fibres, such that the general fibre of $f$ is hyperelliptic of genus 3 and $\si$ induces on it the hyperelliptic involution;
\item the group $\Tors(S)$ contains a subgroup isomorphic to $\Z_2^2$.
\end{enumerate}
\end{prop}
\begin{proof}
Let $C:=K_{W'}+M'$. By Lemma \ref{C}, $|C|$ is a free pencil of rational curves. Notice that $CN'_i=(2K_{W'}+D')N'_i=0$ for every $i$ by Proposition \ref{W'}, so that the curves $N_i'$ are contained in curves of $|C|$.
Since  $C\equiv 2K_{W'}+D'$ and $D'+N_1'+\dots+N'_6$ is divisible by $2$ in $\Pic(W')$ by Proposition \ref{W'}, (iii),   the divisor  $C+N_1'+\dots+N'_6$ is also divisible by $2$. 
Let $Y\to W'$ be the double cover branched on $C+N_1'+\dots+N'_6$, where $C\in |C|$ is general. The surface $Y$ is smooth and the usual formulae for double covers give $\chi(Y)=0$. 
Pulling back $|C|$ to $Y$ one obtains a fibration $h\colon Y\to \Ga$, where $\Ga$ is a smooth curve and the general fibre of $h$ is isomorphic to $\pp^1$. Hence $Y$ is a ruled surface with $q(Y)=1$ and $h$ is the Albanese pencil.

Arguing as in \cite[Theorem 3.2]{dmp}, one shows that there exist effective divisors $A_1,A_2, A_3$ on $W'$ such that, up  to a permutation of the indices, the curves  $2A_1+N'_1+N'_2$, $2A_2+N'_3+N'_4$, $2A_3+N'_5+N'_6$ belong to $|C|$.

We have $CD'=4$, hence by Remark \ref{KC} the system  $|C|$ gives a pencil $|\tilde{C}|$ on $S$ with $K_S\tilde{C}=4$. Since $CN'_i=0$ for every $i$, we have $\tilde{C}^2=0$ and $|\tilde{C}|$ defines a fibration  $f\colon S\to\pp^1$ of hyperelliptic curves of genus 3. The curves of $|C|$ containing the $N'_i$ give rise to double fibres of $f$.

Statement (ii) follows trivially from the existence of three double fibres of $f$.
 \end{proof}
 
 \begin{rem}
In this case it is possible, using the same type of reasoning as in
Corollary 7.6 of \cite{CCM2}, to show that $S$ is a degeneration of surfaces with
non-birational bicanonical map originally
described by Du Val as double planes (cf. \cite{SantaCruz}). Indeed $S$ is  birationally equivalent to a double cover
of $\pp^2$ branched on a curve which is the union of three lines
$r_1,r_2,r_3$ meeting in a point $q_0$ and of a curve of degree 13 with the
following singularities:
\begin{itemize}
\item a 5-uple point at $q_0$;
\item a point $q_i\in r_i$, $i=1,2,3$, of type $[4,4]$, where the tangent
line is $r_i$;
\item three additional  $4-$uple points $q_4,q_5,q_6$ such that there is no conic through $q_1,\dots,q_6$.
\end{itemize}
\end{rem}

\medskip  
 \subsection{The case $K_{W'}^2=-2$} \label{K2W'=-2}

 As in the previous case we consider  $M':=K_{W'}+D'$. 
Recall that ${M'}^2=2$ and $K_{W'}M'= -2$. Moreover $M'$ and $D'$ are nef (cf.  Proposition \ref{W'}).

\begin{lem} \label{G}\begin{enumerate}[\rm(i)]
\item One has $h^0(W', \OO_{W'}(K_{W'}+M'))=1$ and, if $G$ is the unique curve in $|K_{W'}+M'|$, then $GM'=0$.
\medskip

\noindent Moreover, up to a permutation of  the indices $\{1,\dots, 6\}$, one has the following:

\item there are  two possible decompositions of $G$:

a) $G=(2E_1+N'_5)+(2E_2+N'_6)$, where $E_1, E_2$ are $-1-$curves such that $E_1N'_5=E_2N'_6=E_1D'=E_2D'=1$ and the divisors $(2E_1+N'_5)$ and $(2E_2+N'_6)$ are disjoint; or

b) $G=4E_1+3N'_5+2Z_1+N'_6$, where $E_1$ is a $-1-$curve and $Z_1$ is a $-2-$curve 
such that $E_1N'_5=Z_1N'_5=Z_1N'_6=E_1D'=1$
and $E_1Z_1=E_1N'_6=0$;

\item  the divisor $N'_1+\ldots+N'_4$ is even, and it is disjoint from $G$.

\end{enumerate}
\end{lem}

\begin{proof}
We will mimick the proof of Lemma 7.1 in \cite{CCM2}.
The first assertion follows from the long exact sequence
obtained from
$$
0 \to \OO_{W'}(K_{W'}) \to \OO_{W'}(K_{W'}+M') \to \OO_{M'} \to 0
$$
because $W'$ is a rational surface by Proposition \ref{K2W'}, (iii).

By definition of $G$, one has that $G^2=GK_{W'}=-4$ and $GM'=0$.
Therefore, since $M'$ is nef, each component $\theta$ of $G$ is such that $\theta M'=0$ and the intersection
form on the components of $G$ is negative definite.
Since $G^2=-4$, there exists an irreducible curve $E_1$ in $G$
such that $E_1^2<0$ and $E_1G=E_1(K_{W'}+M')<0$.
Since $M'E_1=0$, one has that $E_1K_{W'}<0$, thus $E_1$ is a  $-1-$curve and $E_1G=-1$, $E_1D'=1$.
Recall that $D'\equiv M'-K_{W'}$ is nef, so the irreducible components of $G$ are either $-1-$curves $E$
such that $EG=-1$ and $ED'=1$, or $-2-$curves $Z$ such that $ZG=ZD'=0$.

Since $D'+N'_1+\cdots+N'_6$ is divisible by 2 and $E_1D'=1$,  $E_1$ must  meet one of the $-2-$curves
$N'_i$, say $N'_5$.
Hence $N'_5(G-E_1)=-N'_5E_1<0$, so $N'_5\leq G$ and moreover $E_1N'_5=1$, otherwise we would get $(E_1+N'_5)^2>0$,
a contradiction because the intersection form on the components of $G$ is negative definite.

Similarly $E_1(G-E_1-N'_5)=-1$ implies that $2E_1+N'_5\leq G$.

Recall that $GK_{W'}=-4$, so either $G$ contains another $-1-$curve $E_2$ or $4E_1\leq G$.
Assume the former case.
Then, arguing as before, one sees that $E_2$ meets $N'_i$, for some $i$, and $2E_2+N'_i\leq G$.
If $i=5$, then $(N'_5+E_1+E_2)^2\ge0$, a contradiction.
So we may assume $i=6$.
Finally the negative definiteness implies that $E_1E_2=0$ and
that case a) of statement (ii) occurs, because $(G-2E_1-2E_2-N'_5-N'_6)^2=0$.

Assume now the latter case, i.e.\ $4E_1\leq G$.
Note that $N'_5$ is the only $-2-$curve contained in $G$ that can intersect  $E_1$. Indeed, if $Z\subset G$ is a $-2-$curve such that $E_1Z\geq1$ and $Z\ne N'_5$, then $(2E_1+N'_5+Z)^2\ge0$,  contradicting again the negative definiteness.

Since $E_1G=-1$ and $E_1(4E_1+N'_5)=-3$, one has that $4E_1+3N'_5\leq G$ and 
the components of $G'=G-4E_1-3N'_5$ are $-2-$curves.
Since $N'_5G=0$ and $N'_5 G'=2$, $G$ contains at least a $-2-$curve $Z_1$ with $Z_1N'_5>0$.
Now $N'_5Z_1=1$, otherwise $(N'_5+Z_1)^2\geq0$ gives a contradiction.
Since $Z_1G=0$, we have $2Z_1\leq G'$.
Recall that $Z_1D'=0$ and $D'+N'_1+\ldots+N'_6$ is even, hence $Z_1$ meets another $-2-$curve $N'_i$, say $N'_6$.
Then $N'_6(G-Z_1)=-N'_6Z_1<0$, so $N'_6\leq G'$.
Finally the negative definiteness implies that $N'_6Z_1=1$.
Then we are in  case (ii), b),   because $(G-4E_1-3N'_5-2Z_1-N'_6)^2=0$.

It remains to prove that $N'_1+\cdots+N'_4$ is divisible by 2 in $\Pic(W')$.
Since $D'+N'_1+\cdots+N'_6$ is even, one has that
$2K_{W'}+D'+N'_1+\cdots+N'_6\equiv G+N'_1+\cdots+N'_6$ is also even.
Hence $G+N'_1+\cdots+N'_6\equiv 2(E_1+E_2+N'_5+N'_6)+N'_1+\cdots+N'_4$ is even, in case a),
and $G+N'_1+\cdots+N'_6\equiv 2(2E_1+2N'_5+Z_1+N'_6)+N'_1+\cdots+N'_4$ is even, in case b).
In both cases, one sees that $N'_1+\cdots+N'_4$ is even.

To finish the proof, it is enough to show that $N'_1, \dots, N'_4$  are disjoint from  $E_1$ and $E_2$  in case (ii), a) and from $E_1$ and $Z_1$  in case (ii), b). Arguing as before, this follows easily from the fact that the components of $G$ and the curves $N'_1, \dots, N'_4$ are orthogonal to the nef divisor $M'$.
\end{proof}

By Lemma \ref{G}, there exists a birational morphism $g\colon W'\to X$ such that $X$ is a smooth rational surface, $G$ is the exceptional divisor of $g$ and  $M'=-g^*K_X$.  In particular $-K_X$ is nef and big and $K^2_X=2$. In case (ii), a) of Lemma \ref{G}  the image of $G$ consists of two points $q_5$ and $q_6$ and in case (ii), b)  it is a single point $q$.
\begin{prop}\begin{enumerate}[\rm(i)]
\item There exists a fibration $f\colon S\to\pp^1$ with 2 double fibres, such that the general fibre of $f$ is hyperelliptic of genus 3 and $\si$ induces on it the hyperelliptic involution;
\item The group $\Tors(S)$ contains a subgroup isomorphic to $\Z_2$.
\end{enumerate}
\end{prop}
\begin{proof} For $i=1,\dots, 4$, write $\Delta_i$ for the image of $N'_i$ in $X$. By Lemma \ref{G}, $\De_1+\dots +\De_4$ is again an even set of  disjoint $-2-$curves. By \cite[1.1]{CCM}, there exist a free pencil $|C'|$ of rational curves of $X$ and effective divisors $A_1$, $A_2$ such that, say, $2A_1+\De_1+\De_2$ and $2A_2+\De_3+\De_4$ belong to $|C'|$. The pull back $|C|$ of $C'$ on $W'$ satisfies $CD'=4$, $CN'_i=0$ for $i=1,\dots, 6$, hence it gives a fibration $f\colon S\to \pp^1$ as in statement (i). The curves $2A_1+\De_1+\De_2$ and $2A_2+\De_3+\De_4$ correspond to  two double fibres of $f$.

Statement (ii) follows trivially from the existence of two double fibres of $f$.
\end{proof}

\begin{rem}
As in the previous case,  it is possible, again using the same type of
reasoning
as in Corollary 7.6 of \cite{CCM2}, to show that S is a degeneration of surfaces with
non-birational bicanonical map originally
described by Du Val as double planes (cf. \cite{SantaCruz}). Indeed  $S$ is  birationally equivalent to a double cover
of $\pp^2$ branched on a curve of degree 14  which splits in two distinct
lines $r_1$ and $r_2$
and a curve of degree 12 with the following singularities:
\begin{itemize}
\item the point $q_0=r_1\cap r_2$ of multiplicity 4;
\item a point $q_i\in r_i$, $i=1,2$, of type $[4,4]$, where the tangent
line is $r_i$;
\item two further points $q_3,q_4$ of multiplicity 4
and two points $q_5,q_6$ of type $[3,3]$, such that there is no conic
through $q_1,\ldots,q_6$.
\end{itemize}

The point $q_6$ is infinitely near to $q_5$, in case  (ii), b)  of Lemma
\ref{G}.
\end{rem}

\medskip
 
 \subsection{The case $K_{W'}^2=-3$} \label{K2W'=-3}

Denote $M':=K_{W'}+D'$. We have ${M'}^2=1$, $K _{W'}M'=-3$
 and $|M'|$ is 2-dimensional. Since $|M'|$
 has no fixed part by Lemma \ref{fixed}, the map
 $\phi_{M'}\colon W'\to \pp^2$ is a birational morphism.  It is an
 easy exercise seeing that the branch curve
 is mapped to a plane curve of  degree 10, which,
as it is well known, has  6 singular points  of type $[3,3]$
(possibly infinitely near).

The original construction proposed by Campedelli (\cite{Cam}) is one of 
these surfaces. For a discussion of possible branch loci and relations 
with the $2$-torsion of $S$ see \cite{St} and \cite{werner}.

 \subsection{The case $K_{W'}^2=-4$} \label{K2W'=-4}

\medskip

We start by noticing that in this case $W'=W$, because $-4\leq K_W^2$ by Proposition \ref{properties}, (iii).

Denote $M:=K_{W}+D$. Recall that  $|M|$ has no fixed part, by Lemma \ref{fixed}. Then  $M^2=0$ and $h^0(W, M)=3$ imply that 
 $|M|=|2C|$, where $|C|$ is a pencil without base points. Since  $K_W^2=-4$ and $K_WD=0$,  we have $K_WC=-2$, hence $|C|$ is a pencil of rational curves. 
Since $CN_i=0, i=1,...,6$ and $DC=2$, $C$
 gives rise to a genus 2 fibration $\tilde{C}$ on $S$ such that $\si$ restricts to the hyperelliptic involution on the general $C$.

Notice that in this case the curve $B_0$ on $W$ must be reducible, because
by Proposition \ref{godeaux} $p_a(R)=-1$ and, of course, $p_a(B_0)=p_a(R)$.
In fact, recalling that $D=2K_W+B_0$, we obtain $B_0^2=-12$, $K_WB_0=8$, hence $p_a(B_0)=-1$.

\begin{rem}
Conversely, assume that the numerical Campedelli surface $S$ has a free pencil $|C|$ of curves of genus 2 and let $\si$ be the involution of $S$ that induces the hyperelliptic involution on the general $C$. Then the results of \cite[\S 1, 2]{xiao} (cf. Remark 2.4, ibidem) show that we have $K^2_{W'}=-4$ in this case.
\end{rem}

\begin{rem}
In this case by  \cite[\S 2]{xiao}  the relative canonical map of $S$
expresses $S$ as a double cover of $\F_0=\pp^1\times \pp^1$ branched in a
curve of degree $(6,8)$, which in the general case has $6$ distinct
singular points of  type $[3,3]$.
\end{rem}
 
 \section{Involutions not composed with the bicanonical map}
 
 In this section we consider case II) of Proposition \ref{cases}, namely here we assume that $k=4$ and  the bicanonical map $\fie\colon S\to\pp^2$ is not  composed with $\si$. 
 We recall that by Proposition \ref{properties} in this case we have $D^2=4$, $K_WD=2$, $-2\le K^2_W\le 1$.

\begin{lem}\label{compB}   Set $m=1-K^2_W$. Then $B_0=\Ga_1+\dots +\Ga_m$,
where the $\Ga_i$ are disjoint $-4-$curves.

\end{lem}
\begin{proof} Notice first of all that $B_0D=D^2-2DK_W=0$.  Let $\Ga$ be an irreducible component of $B_0$ and write $p^*\Ga=2\tilde{\Ga}$. We have $D\Ga=0$, since $D$ is nef, and thus $\epsi^*K_S\tilde{\Ga}=0$, since $p^*D=\epsi^*(2K_S)$ by Proposition \ref{properties}. Since $\tilde{\Ga}$ is disjoint from the exceptional locus of $\epsi$ by construction, it follows that $\tilde{\Ga}$ is a $-2-$curve. Hence $\Ga^2=-4$  and $\Ga$ is a smooth rational curve. 
Now let $m\ge 0$ denote the number of components of $B_0$. By the adjunction formula we have $K_WB_0=2m$. On the other hand, we can compute:
$$2m=K_WB_0=K_W(D-2K_W)=2-2K^2_W.$$
Finally, the components of $B_0$ are disjoint, since $B_0$ is smooth. 
\end{proof}
\begin{cor} If $K^2_{W}\le 0$, then $K_S$ is not ample.
\end{cor}
\begin{proof} By Lemma \ref{compB}, the branch divisor $B$ of the map $\pi\colon S\to\Si$ contains at least  a smooth rational curve $\Ga$ with $\Ga^2=-4$. Then the inverse image of $\Ga$ in $S$ is a $-2-$curve and $K_S$ is not ample.
\end{proof}

\begin{prop} \label{W} We have the following possibilities:
\begin{enumerate}[\rm(i)]
\item $K^2_W=1$, $W$ is minimal of general type and $B_0=0$;
\item $K^2_W=0$, $W$ is minimal and properly elliptic;
\item $K^2_W=-1, -2$ and $W$ is not of general type.
\end{enumerate}
\end{prop}
\begin{proof} Recall that $-2\le K^2_W\le 1$ by Proposition \ref{properties}. 
If $K^2_W=1$, then by Lemma \ref{compB} we have $B_0=0$ and $K_S=\pi^*K_{\Si}$, hence $K_W=\eta^*K_{\Si}$ is nef and big and $W$ is  minimal  of general type. 

Next we show that if $K^2_W\le 0$, then $W$ is not of general type. So assume by contradiction  that $W$ is of general type.  Let $t\colon W\to W_1$ be the morphism to the minimal model and write $K_W=t^*K_{W_1}+E$, where $E>0$. 
Since $DK_W=2$ and $D$ is nef, we have $Dt^*K_{W_1}\le 2$. On the other hand, since $K_{W_1}^2>0$ the index theorem applied to $D$ and $t^*K_{W_1}$ gives $Dt^*K_{W_1}\ge 2$. So we get $Dt^*K_{W_1}=2$ and $D\sim 2t^*K_{W_1}$. This implies $B_0+2E\sim 0$, a contradiction since $B_0+2E >0$.

Assume now that $K^2_W=0$. By Lemma \ref{compB}, $B_0$ is a smooth rational curve with $B_0^2=-4$.
By the exact sequence:
\begin{equation}
0\to H^0(2K_W)\to H^0(2K_W+B_0)\to H^0(\OO_{B_0})
\end{equation}
we obtain $1\le h^0(2K_W)\le 2$, hence $W$ has nonnegative Kodaira dimension.
We have seen that $W$ is not of general type, hence it is minimal and it is either properly elliptic or Enriques. Since $K_WD=2\ne0$, the latter case does not occur. This finishes the proof.
\end{proof}

\begin{rem}
 By Proposition \ref{W},  the desingularization $W$ of the quotient surface 
$S/\si$
may  be a numerical Godeaux 
surface, an elliptic surface, birational to an Enriques surfaces or 
rational.

Unlike the previous case, in which one knows examples for all the 
possibilities for $W$,  in this case  we do not know any example for which 
$W$ is rational. R. Barlow in \cite{Ba1}, \cite{Ba2} presents examples of numerical Godeaux surfaces with four nodes double covered by numerical Campedelli surfaces and  the new examples we present  such that $W$ is not a 
surface of 
general type  are  obtained by specializing one of these constructions (cf. \S 5).

It is possible to make  a more detailed analysis of the cases with 
$K^2_W\le 0$, in the style of the previous section, but since the arguments
 are very lengthy  and all the examples we know are obtained by specialization, we do not think 
worthwhile
including it here.

\end{rem}

\section{Examples}

In this section we study   some families of numerical Campedelli surfaces with an involution,  providing examples for  the  cases 3.2 to 3.5 in \S 3 and for the cases (i)---(iii)  in Proposition \ref{W}. 
\bigskip

{\bf Example 1.} {\em Numerical Campedelli surfaces with torsion $\Z_2^3$.}

These surfaces have two different descriptions: as the quotient by a free $\Z_2^3-$action of the intersection of four quadrics in $\pp^6$ (cf. \cite{miyaoka1}, \cite{reid})  and as $\Z_2^3-$covers of $\pp^2$ branched on 7 lines (cf. \cite{kulikov}). We use the second description, which is more suitable for our purposes.
Two special instances of surfaces in this family are the Burniat surface with $K^2=2$ and the classical Campedelli surface (cf. \cite[\S 4]{kulikov}).

Set $G:=\Z_2^3$ and let $\chi_1,\chi_2, \chi_3$ be generators of $G^*$, the group of characters  of $G$. By \cite[Proposition 2.1 and Corollary 3.1]{ritaabel}, to give a normal $G-$cover  $p\colon X\to \pp^2$ it is enough to give an effective divisor $D_g$ for every $0\ne g\in G$ and line bundles $L_1, L_2, L_3$ on $\pp^2$ such that the divisor $\De:=\sum_{g\ne 0} D_g$ is reduced and the following relations are satisfied:
$$ 2L_i\equiv \sum_{g\ne 0}\epsi_i(g)D_g, \quad i=1,2,3,$$
where we define $\epsi_i(g)=0$ if $\chi_i(g)=1$ and  $\epsi_i(g)=1$ if $\chi_i(g)=-1$.

Here we take the $D_g$ to be distinct lines in $\pp^2$ and we set $L_i:=\OO_{\pp^2}(2)$, $i=1,2,3$.
Moreover we make the following assumptions on the configuration  of the lines $D_g$:

1) at most three of the $D_g$ pass through the same point;

2) if $D_{g_1}$, $D_{g_2}$, $D_{g_3}$ pass through the same point, then $g_1+g_2+g_3\ne 0$.

We now examine the singularities of $X$.   By \cite[Proposition 3.1]{ritaabel}, $X$ is singular above a point $P\in\pp^2$ if and only if $P$ lies on three branch lines $D_{g_1}$, $D_{g_2}$ and $D_{g_3}$. To resolve the singularity, let $\psi\colon \hat{\pp}\to \pp^2$ be the blow up of $\pp^2$ at $P$, let $E$ be the exceptional curve of $\psi$ and consider the $G-$cover $\hat{p}\colon \hat{X}\to\hat{\pp}$ obtained from $p$ by base change and normalization. Write $g_0:=g_1+g_2+g_3$. By \cite[\S 3]{ritaabel}, the components of the branch divisor of $\hat{p}$ are the following:
$\hat{D}_g:=\psi^*D_g$ if $g\ne g_0,\dots, g_3$, $\hat{D}_{g_0}:=\psi^*D_{g_0}+E$, $\hat{D}_{g_i}:=\psi^*D_{g_i}-E$ for $i=1,2,3$. 
The surface $\hat{X}$ is smooth above $E$ and $\psi\inv(E)$ is a $-2-$curve.
Hence $X$ has a rational double point of type $A_1$ over $P$.
 We have $2K_X=\psi^*(\OO_{\pp^2}(1))$, hence $K_X$ is ample and $X$ is the canonical model of a surface $S$ of general type  with $K^2_S=2$. By the projection formulae for abelian covers we have $|2K_X|=\psi^*|\OO_{\pp^2}(1)|$, hence $h^0(X,2K_X)=3$, $\chi(S)=\chi(X)=1$ and $\psi$ is the bicanonical map of $X$.

Kulikov (\cite[Thm. 4.2]{kulikov}) shows that the automorphism group of the general surface in this family coincides with the Galois group $G=\Z_2^3$ of the bicanonical map.
The result that follows is a partial refinement of his, and gives evidence for the difficulty of finding an involution of a numerical Campedelli surface such that
its bicanonical map is not composed with it.
\begin{prop}\label{invo}
Let $S$ be a numerical Campedelli surface with torsion $\Z_2^3$ and let $\si$ be an involution of $S$.
Then $\si$ is in the Galois group $G=\Z_2^3$ of the bicanonical map of $S$.
\end{prop}
\begin{proof}
Assume by contradiction that $\si$ is an involution of $S$ such that the bicanonical map $\fie\colon S\to\pp^2$ is not composed with $\si$. Since $G$ is defined intrinsically, we have $\si G\si=G$ and $\si$ induces an involution of $\pp^2$ that we denote by   $\bar{\si}$. 
Since the set of lines $D_g$ contains at least 4 lines in general position, $\bar{\si}$ induces a non trivial permutation of the $D_g$. Denote by $h$ the automorphism of $G$ defined by $h(g)=\si g\si$. Then we have $\si(D_g)=D_{h(g)}$, and 
it follows that $h$ is a non trivial automorphism of $G$. Since $h$ has order 2, we can find generators $e_1, e_2,e_3$ of $G$ such that $h(e_i)=e_i$ for $i=1,2$ and $h(e_3)=e_3+e_1$. Hence the lines $D_{e_1}, D_{e_2}, D_{e_1+e_2}$ are fixed for $\bar{\si}$, while $D_{e_3}$ and $D_{e_3+e_1}$ are exchanged by $\bar{\si}$ and the same happens to $D_{e_3+e_2}$ and $D_{e_3+e_2+e_1}$.
Then, taking also  into account the combinatorial conditions on the configuration of the lines $D_g$, up to   exchanging $e_2$ and $e_1+e_2$, we can find  homogeneous coordinates on $\pp^2$ such that $\bar{\si}(x_0,x_1,x_2)=(x_0, x_1,-x_2)$ and such that:
\begin{gather*} D_{e_1}=\{x_1=0\}, \ D_{e_2}=\{x_0=0\}, \ D_{e_1+e_2}=\{x_2=0\}\\
D_{e_3}=\{ax_0+bx_1+cx_2=0\},\  D_{e_3+e_1}=\{ax_0+bx_1-cx_2=0\},\\
D_{e_3+e_2}=\{a'x_0+b'x_1+c'x_2=0\}, \ D_{e_3+e_1+e_2}=\{a'x_0+b'x_1-c'x_2=0\}.
\end{gather*}
Since  $h$ maps  the subgroup $H$ of $G$ generated by $e_1$ and $e_3$ to itself, $\si$ induces an involution  of the surface $Z:=S/H$ that lifts $\bar{\si}$. 
On the other hand, the function field $\C(Z)$ is the quadratic extension of $\C(\pp^2)$ obtained by adding the square root of $x_0(a'x_0+b'x_1+c'x_2)(a'x_0+b'x_1-c'x_2)/x_2^3$, and it is easy to check that the action of $\bar{\si}$ on $\C(\pp^2)$ cannot be extended to an automorphism of order 2 of $\C(Z)$. Hence we have obtained a contradiction.
\end{proof}
We  now study the involutions of $G$. There are different cases, according to the relative positions of the lines in $\De$.
\medskip

\underline{Case 1:}  the lines of $\De$ are in general position.

In this case $K_S$ is ample and therefore $W=W'$ for any involution of $S$ by Remark \ref{remW'}.

The divisorial part $R$ of the fixed locus on $S$ of any $0\ne \si\in G$  is a paracanonical curve. Hence,  the adjunction formula $K_S\equiv \pi^*K_{\Si}+R$ gives  that $2K_{\Si}\equiv 0$ and  $\Si$  is an Enriques surface with 6 nodes. So this is an instance of case 3.1 of \S 3.
Other examples of this case, with torsion $\Z_2\times \Z_4$, appear in \cite{Naie}.
\medskip 

\underline{Case 2:}  the divisor $\De$ has one triple point $P$, lying on the lines $D_{g_1}$, $D_{g_2}$, $D_{g_3}$. Consider the involution $g_0:=g_1+g_2+g_3$. In this case the cover $\hat{p}\colon \hat{X}\to\hat{\pp}$ is smooth and  we have $S=\hat{X}$. The divisorial part of the fixed locus of $g_0$ on $S$ is the disjoint union of the $-2-$curve that resolves the singularity of $X$ and of a paracanonical curve.  Hence one gets $K^2_W=-1$. Since the only $-2-$curve of $S$ is in the fixed locus of $g_0$, we have $W'=W$ and the surface $W$ is rational by Proposition \ref{K2W'}, namely this is an example of case 3.2.
In fact, it is easy to check that the lines through the point $P\in \pp^2$ pull back to a pencil of rational curves on $W$, which in turn gives a free pencil of hyperelliptic curves of genus 3 with three double fibres on $S$.
\medskip

\underline{Case 3:} the divisor $\De$ has a triple point $P_1$, lying on the lines $D_{g_1}$, $D_{g_2}$, $D_{g_3}$, and another triple point $P_2$, lying on the lines $D_{h_1}$, $D_{h_2}$ and $D_{h_3}$, with $g_1+g_2+g_3=h_1+h_2+h_3=:g_0$.

Arguing as in Case 2, one shows that the fixed locus of $g_0$  on $S$ is the disjoint union of a paracanonical curve and of the two $-2-$curves that resolve the double points of $X$ lying above $P_1$ and $P_2$. We have $W=W'$, $K^2_W=-2$ and  $W$ is rational. Hence this is an example of case 3.3.
\medskip

\underline{Case 4:} the divisor $\De$ has three  triple points: $P_1$, lying on the lines $D_{g_1}$, $D_{g_2}$, $D_{g_3}$,  $P_2$, lying on the lines $D_{h_1}$, $D_{h_2}$ and $D_{h_3}$, and $P_3$, lying on the lines $D_{f_1}$, $D_{f_2}$ and $D_{f_3}$. Moreover, we assume that  $g_1+g_2+g_3=h_1+h_2+h_3=f_1+f_2+f_3=:g_0$. We remark that  the existence of such a configuration of lines is not difficult to verify.

Arguing as in Case 2, one shows that the fixed locus of $g_0$  on $S$ is the disjoint union of a paracanonical curve and of the three $-2-$curves that resolve the double points of $X$ lying above $P_1$, $P_2$ and $P_3$. We have $W=W'$, $K^2_W=-3$ and  $W$ is rational. Hence this is an example of case 3.4.

\begin{rem}One can check that  $\De$ cannot have  four triple points $P_1,\dots, P_4$ such that $P_i$ lies on $D_{g_{1i}}$, $D_{g_{2i}}$, $D_{g_{3i}}$ with $g_{1i}+g_{2i}+g_{3i}=g_{1j}+g_{2j}+g_{3j}\ne 0$ for every $i, j=1,\dots, 4$. Hence, by Proposition \ref{invo} the cases 1--4  described above are essentially the only possibilities for an  involution of a numerical Campedelli  surface with torsion $\Z_2^3$. 
\end{rem}

{\bf Example 2.} {\em A family of  numerical Campedelli surfaces with torsion $\Z_3^2$ and two involutions.}

This example has been kindly communicated to us by JongHae Keum, who attributes it to Xiao Gang and Beauville.

Consider $X:=\pp^2\times \pp^2$ with homogeneous coordinates $(x_0,x_1,x_2; y_0,y_1,y_2)$ and let two generators $g_1$ and $g_2$ of the group $G:=\Z_3^2$ act on $X$ as follows:
\begin{gather*}
(x_0,x_1,x_2; y_0,y_1,y_2)\overset{g_1}{\mapsto}(x_1,x_2,x_0; y_1,y_2,y_0);\\(x_0,x_1,x_2; y_0,y_1,y_2)\overset{g_2}{\mapsto}(x_0,\omega x_1, \omega^2x_2; y_0,\omega^2 y_1,\omega y_2),
\end{gather*}
where $\omega$ is a primitive $ 3-$rd root of $1.$
Consider the family of surfaces $Y$ of $X$ defined by the equations:
\begin{eqnarray*}
x_0y_0+x_1y_1+x_2y_2=0,\\ (x_0^3+x_1^3+x_2^3)(y_0^3+y_1^3+y_2^3)+\lambda x_0x_1x_2y_0y_1y_2=0.
\end{eqnarray*}
For a general value  of the parameter $\lambda\in \C$ the surface $Y$ is smooth and simply connected with $K^2_Y=18$, $p_g(Y)=8$, and the group $G$ acts freely on it. Hence the quotient surface $S:=Y/G$ is a numerical Campedelli surface with fundamental group equal to $G$. 

The surface $Y$ is mapped to  itself also by the involution $\tilde{\si}_1$ of $X$ defined by:
$$(x_0,x_1,x_2; y_0,y_1,y_2)\overset{\tilde{\si}_1}{\mapsto} (y_0,y_1,y_2; x_0,x_1,x_2).$$
The involution $\tilde{\si}_1$ satisfies the following relations:
\begin{equation}\label{rel}
\tilde{\si}_1 g_1=g_1\tilde{\si}_1, \ \tilde{\si}_1g_2=g_2^2\tilde{\si}_1,
\end{equation}
hence $G$ and $\tilde{\si}_1$ generate a group $G_0$ of order 18, the involution $\tilde{\si}_1$ induces an involution $\si_1$ of $S$ and we have  $Y/G_0=S/\si_1$.

The fixed locus of $\tilde{\si}_1$ on $Y$ consists of $12$ points and the same is true for $\tilde{\si}_1g_2$ and $\tilde{\si}_1g_2^2$, since these involutions are conjugated to $\tilde{\si}_1$. Consider now an element of $G_0$ of the form $\tilde{\si}_1g$, where $g\in G\setminus\langle g_2 \rangle$. The relations (\ref{rel}) imply that $(\tilde{\si_1}g)^2$  is a nonzero element of $G$, hence in particular $\tilde{\si}_1g$ has no fixed points on $Y$.
It follows that $\si_1$ has 4 fixed points on $S$ and the quotient surface $T:=S/\si_1$ is a numerical Godeaux surface. By \cite[\S 0]{Ba1}, the fundamental group of $T$ is isomorphic to $\Z_3$.
 Hence we have an example in which the bicanonical map is not composed with the involution and  the quotient surface is of general type, that is case (i) of Proposition \ref{W}.
\begin{comment}
 \begin{rem} Similar examples are contained in \cite{Ba2}, where the quotient surfaces have fundamental group equal to $\Z_2$ or to $\Z_4$. Moreover, the famous Barlow surface, which is simply connected, can also be obtained as the quotient of a numerical Campedelli surface with torsion $\Z_5$ by an involution (cf. \cite{Ba1}).
\end{rem}
\end{comment}
 
 Consider now the involution $\tilde{\si}_2$ defined by:
$$(x_0,x_1, x_2;y_0,y_1,y_2)\overset{\tilde{\si}_2}{\mapsto} (x_0, x_2,x_1; y_0,y_2,y_1).$$
For every $g\in G$ one has the relation $\tilde{\si}_2 g=g\inv {\tilde{\si_2}}$. Hence the group $G_0$ generated by $G$ and $\tilde{\si}_2$ has order 18. $G_0$ contains nine elements of order 2, that form a conjugacy class.
The surface $Y$ is mapped to itself by $\tilde{\si}_2$ and $\tilde{\si}_2$ induces an involution of $S$ that we denote by $\si_2$. We have $Y/G_0=S/\si_2$.

The fixed locus of $\tilde{\si}_2$ on the threefold $\{x_0y_0+x_1y_1+x_2y_2=0\}\subset X$ consists of 3 disjoint rational  curves:
$\Ga_1=\{(0,1, -1; a, b,b)| (a,b)\in\pp^1\}$, $\Ga_2=\{( a, b,b; 0,1,-1)| (a,b)\in\pp^1\}$, $\Ga_3=\{(a,b, b; -2b, a,a)| (a,b)\in\pp^1\}$. 
It is not difficult to check that $\Ga_1$ and $\Ga_2$ are contained in $Y$, while $\Ga_3$ meets the general $Y$ at 6 distinct points. Since $K_Y$ is the restriction of $\OO_{\pp^2\times\pp^2}(1,1)$ to $Y$, we have $K_Y\Ga_i=1$, for $i=1,2$ and $\Ga_1$, $\Ga_2$ are $-3-$curves on $Y$.
Hence the fixed locus of $\si_2$ on $S$ is the union of 6 isolated points and two $-3-$curves, and thus $K^2_W=-4$. If $Y$ is smooth, then $K_Y$ and $K_S$ are ample and we have $W=W'$ by Remark \ref{remW'}. So this is also  an example of case 3.5.

We are now going to show that the involution $\si_2$ of $S$ is actually induced by a genus 2 pencil, as explained in 3.5. 
Consider the pencil of hypersurfaces of $X$ spanned by $x_0x_1x_2$ and $x_0^3+x_1^3+x_2^3$ and denote by $|F|$ the restriction of this pencil to $Y$. The fixed part of $|F|$ is the union of the curves in  the orbit of $\Ga_1$ under the group $G$. Then we can write $|F|=Z+|C|$, where the $Z$ is the  disjoint union of nine $-3-$curves and $|C|$ has no fixed part. On the surface $Y$ we have $F^2=27$, $K_YF=27$.   Using this and $F\Ga=0$ for every component $\Ga$ of $Z$,  one gets $C\Ga=3$, $K_YC=18$, $C^2=0$. Every element of $|C|$ is mapped to itself by  $G_0$. Hence $|C|$ induces a genus 2 pencil $|C'|$ of $S$  such that every element of $|C'|$ is mapped to itself by ${\si}_2$. Finally, using $\Ga_2F=3$ and the fact that the $G-$orbits of $\Ga_1$ and $\Ga_2$ are disjoint, we get $C\Ga_1=C\Ga_2=3$. So the general $C'$ meets the fixed locus of ${\si}_2$ at 6 points, hence ${\si}_2$ restricts to the hyperelliptic involution on $C'$.
\medskip

{\bf Example 3.}  {\em Numerical Campedelli surfaces with an involution with which the bicanonical map is not composed and such that the quotient is not of general type.}

Here we provide examples for cases (ii) and (iii) of Proposition \ref{W}.
These examples are obtained by specializing a construction due to R. Barlow (\cite{Ba2}). We start by recalling  briefly her construction.

Consider the space $\pp^6$ with homogeneous coordinates $(x_1,\dots, x_7)$ and the automorphisms of $\pp^6$ defined as follows:
\begin{gather*}
(x_1,\dots, x_7)\overset{t}{\mapsto}(\zeta x_1, \zeta^2x_2,\dots, \zeta^7x_7)\\
(x_1,\dots , x_7)\overset{a}{\mapsto} (x_3, x_6, x_1, x_4, x_7, x_2, x_5),
\end{gather*}
where $\zeta$ is a primitive $8-$th root of 1. The automorphism $t$ has order $8$, the automorphism $a$ has order 2 and one has: $$ata=t^3.$$
 Hence $a$ and  $t$ generate a subgroup $G$ of  order 16 of $\Aut(\pp^6)$.
Consider the intersection $Y$ of the following four quadrics of $\pp^6$:
\begin{gather*}
F_0:=b(x_1x_7+x_3x_5)+ax_4^2+fx_2x_6\\
F_2:=cx_1^2+dx_3x_7+ex_4x_6+hx_5^2\\
F_4:=k(x_2^2+x_6^2)+gx_1x_3+mx_5x_7\\
F_6:=cx_3^2+dx_1x_5+ex_4x_2+hx_7^2
\end{gather*}
Barlow proves that for a general  choice of the coefficients $a$, $b$, $c$, $d$, $e$, $f$, $g$, $h$, $k$, $m$ the following are true:
\begin{itemize}
\item $Y$ is a smooth surface mapped to itself by $G$;
\item the subgroup $\Z_8$ of $G$ generated by $t$ acts freely on $Y$;
\item the involution $a$ has 8 isolated fixed points on $Y$.
\end{itemize}
It follows easily from the properties above that the quotient surface $S:=Y/\Z_8$ is a numerical Campedelli surface with torsion $\Z_8$. In addition, the involution $a$ of $Y$ induces an involution $\si$ of $S$ with four isolated fixed points. The quotient surface $\Si:=S/\si$ has four nodes and its minimal desingularization $W$ is a minimal surface of general type with $K^2_W=1$ and $p_g(W)=0$, namely a numerical Godeaux surface. Barlow also shows that $\pi_1(W)=\Z_2$.

Let  $\Ga$ be the group of automorphisms of $\pp^6$ of the form $\mbox{Diag}(1, \lambda, 1, \mu, \nu, \lambda,\nu)$ for $\lambda, \mu,\nu\in\C^*$.   The elements of $\Ga$ commute with $a$ and $t$ and act on the family of surfaces $Y$, hence the family of numerical Campedelli surfaces $S$ that we obtain has at most 4 moduli.

We are going to specialize this construction by letting $S$ acquire one or two ordinary double  points which are fixed by $\si$ and whose images in $\Si$ are quotient singularities of type $\frac{1}{4}(1,1)$.
Passing to the minimal desingularization $S'$ of $S$ we obtain an involution whose fixed locus consists of four isolated points and of the $-2-$curves that resolve the singularities of $Y$. In the case of one double point we get an example of case (ii) of Proposition \ref{W}. In particular, the minimal desingularization $W$ of $\Si$ is a properly elliptic surface. In the case of two double points we have an example of case (iii) of Proposition \ref{W} and we will show that $W$ is a non minimal Enriques surface.

The fixed locus of $a$ on $\pp^6$ consists of the $\pp^3$ defined by $$x_1-x_3=x_2-x_6=x_5-x_7=0$$ and of the $\pp^2$ defined by $$x_1+x_3=x_2+x_6=x_5+x_7=x_4=0.$$ In \cite{Ba2} it is shown that the general $Y$ intersects the $\pp^3$ in 8 points and it does not intersect the $\pp^2$.
Let $P_1\in \pp^2$ be the point $(1,1, -1,0, 1,-1,-1)$ and let $P_2:=t^4(P_1)\in\pp^2$ the point $(1,-1,-1,0,1,1,-1)$. Let $P_3, \dots, P_8$ denote the remaining points in the orbit of $P_1$ under the action of $\Z_8$. 
The surfaces  $Y$ that contain $P_1,\dots, P_8$ are  defined by  four quadrics as follows:
\begin{gather*}
F_0:=b(x_1x_7+x_3x_5)+ax_4^2-2bx_2x_6\\
F_2:=cx_1^2+dx_3x_7+ex_4x_6-(c+d)x_5^2\\
F_4:=k(x_2^2+x_6^2)+gx_1x_3+(2k-g)x_5x_7\\
F_6:=cx_3^2+dx_1x_5+ex_4x_2-(c+d)x_7^2
\end{gather*}
It is easy to verify that the tangent space to the general $Y$ at $P_1$ has dimension 3 and $a$ acts on it as multiplication by $-1$. Since the points $P_1,\dots, P_8$ form an orbit under the action of $t$ and $Y$ is mapped to itself by $t$, the singularities of $Y$ at $P_1, \dots, P_8$ are isomorphic.
\begin{rem}\label{sing}
 The orbit of $P_1$ under the action of $\Ga$ is dense in the $\pp^2$ fixed by $a$. It follows that if a surface $Y$ intersects this $\pp^2$ in a point $P$, then $Y$ is singular at $P$ and the subspace of the  tangent space to $Y$ at $P$ on which $a$ acts as multiplication by $-1$ has dimension at least $3$. In addition, if $P$ satisfies $x_1x_2x_5\ne 0$ and $Y$ is general among the surfaces through $P$, then  the tangent space to $Y$ at $P$ has dimension 3 and $a$ acts on it as multiplication by $-1$.
\end{rem}

We claim that for a general choice of the parameters $a$, $b$, $c$, $d$, $e$, $g$, $k$ the surface $Y$ satisfies the following conditions:

1) the subgroup generated by $t$ acts freely on $Y$;

2) $Y$ meets the $\pp^3$ fixed by $a$ in 8 points and it meets the $\pp^2$ fixed by $a$ in $P_1$ and $P_2$;

3) $Y$ has an ordinary double point in $P_1,\dots, P_8$ and it is smooth elsewhere. 

Conditions 1)--3) are open, hence it is enough to check them for  one surface $Y$.
Let $Y_0$ be the surface corresponding to the following choice of parameters:
$$a=e=-1, b= c=d=g=k=1.$$

Using a computer program (we have used Singular) one checks the following:
\begin{itemize}
\item $Y_0$ does not intersect the spaces $H_1:=\{x_1=x_3=x_5=x_7=0\}$ and $H_2:=\{x_2=x_4=x_6=0\}$ fixed by $t^4$, hence condition 1) is satisfied;
\item $Y_0$ intersects the $\pp^3$ fixed by $Y$ at 8 points;
\item  the scheme of singular points of $Y_0$ has dimension 0 and degree 8.
\end{itemize}

Since we already know that $Y$ is singular at  $P_1,\dots , P_8$, the last condition above implies 3).  The fact that $Y_0$ meets the $\pp^2$ fixed by $a$ only at $P_1, P_2$ is now a consequence of Rem. \ref{sing}. 
Hence conditions 1)--3) are satisfied by $Y_0$ and therefore they are satisfied by the general $Y$ that has nonempty intersection with the $\pp^2$ fixed by $a$.  For such a surface $Y$, the quotient surface $S:=Y/\Z_8$ has an ordinary double point at the image point $P$ of $P_1,\dots,P_8$ and it is smooth elsewhere. Hence $S$ is the canonical model of a numerical Campedelli surface.
Let $S'$ be the minimal resolution of $S$, let $Z$ be the exceptional curve  and let $\si'$ be the involution of $S'$ induced by $a$.
Since $a$ acts on the tangent space to $Y$ at $P_1$ as multiplication by $-1$, the fixed locus of $\si'$ on $S'$ consists of the curve $Z$ and of 4 isolated fixed points. Hence we have $K^2_W=0$ by  Lemma \ref{compB} and $W$ is minimal and  properly elliptic by Proposition \ref{W}. 
Applying the argument used in \cite{Ba2}, one can show that the fundamental group of $W$ is $\Z_2$.

  Since the elements of $\Ga$ with $\lambda=\nu=1$ act on the family of surfaces $Y$ passing through $P_1$, the family of Campedelli surfaces with one node that we have constructed has at most 3 moduli.

We are now going  to degenerate the construction further, letting $S$ acquire two double points,  and thus obtain an example with $K^2_W=-1$. 
\begin{comment}Now we choose $P_1=(1,1,-1, 0,-1,-1, 1)$, $P_2=t^4P_1=(1,-1,-1,0, -1, 1, 1)$,  $Q_1:=(1, 2, -1, 0, -4, -2, 4)$, $Q_2:=t^4Q_1=(1, -2, -1, 0, -4, 2, 4)$. We denote  by $P_3, \dots , P_8$ the remaining points in the orbit of $P_1$ under the action of $t$ and  by $Q_3, \dots , Q_8$ the remaining points in the orbit of $Q_1$. Notice that the points $P_1$, $P_2$, $Q_1$ and $Q_2$ lie in the $\pp^2$ fixed by $a$. 
\end{comment}
Set $Q_1:=(1,2,-1,0,4,-2,-4)$ and $Q_2:=t^4Q_1=(1,-2,-1,0,4, 2,-4)$ and denote by $Q_3,\dots,Q_8$ the remaining points in the orbit of $Q_1$ under the action of $t$. 
The surfaces $Y$ through $P_1,\dots ,P_8$ and $Q_1,\dots ,Q_8$ are defined by the following   four quadrics:
\begin{gather*}
F_0:=b(x_1x_7+x_3x_5)+ax_4^2-2bx_2x_6\\
F_2:=cx_1^2-\frac{5}{4}cx_3x_7+ex_4x_6+\frac{1}{4}cx_5^2\\
F_4:=5(x_2^2+x_6^2)+8x_1x_3+2x_5x_7\\
F_6:=cx_3^2-\frac{5}{4}cx_1x_5+ex_4x_2+\frac{1}{4}cx_7^2
\end{gather*}

By Rem. \ref{sing} every surface as above is singular at $P_1, \dots, P_8, Q_1, \dots, Q_8$.
Let now $Y_0$ be the surface corresponding to the following choice of parameters:
$$a=-1, b=1, c=4, e=-1.$$
Also in this case, we have used the computer program  Singular to check that $Y_0$ has the following properties:
\begin{itemize}
\item the automorphism $t$ acts freely on $Y_0$;
\item $Y_0$ intersects the $\pp^3$ fixed by $a$ at 8 points;
\item $Y_0$ intersects the $\pp^2$ fixed by $a$ at $Q_1,Q_2, P_1,P_2$;
\item the scheme of singular points of $Y_0$ has dimension 0 and degree 16, and thus $Y_0$ has  a node at $P_1, \dots, P_8, Q_1,\dots , Q_8$ and is smooth elsewhere.
\end{itemize}
Since these properties are open, they hold for the general $Y$ passing through $P_1$ and $Q_1$. 
The quotient surface $S:=Y/\Z_8$ has two nodes which are fixed by $a$ and $K_S$ is ample. Let $S'$ be the minimal  desingularization of $S$, let $Z_1$ and $Z_2$ be the exceptional curves on $S'$ and let $\si$ be the involution of $S$ induced by $a$. The fixed locus of $\si$ on $S$ consists of 4 isolated points and of the curves $Z_1$ and $Z_2$ (cf. Remark \ref{sing}). Hence we have $K^2_W=-1$ and this is an example of case (iii) of Proposition \ref{W}.

As in the case of one node, one can use the same argument as in \cite{Ba2} to show that $\pi_1(W)=\Z_2$. Hence $W$ is not rational and, by Proposition \ref{W},  it is birational either to an  Enriques surface or to a properly elliptic surface. 
We are going to see that in fact $W$ is birational to an Enriques surface.

The intersection of $Y$ with the hypersurface $x_4^2=0$ is a bicanonical curve  which descends to a bicanonical curve $2C\subset S$ passing through the nodes of $S$. Pulling back to $S'$ we obtain a bicanonical curve $2C'=2Z_1+2Z_2+2G'$, where $G'$ is effective. 

 By the adjunction formula, there is an effective divisor $G$ on $W$ such that $G\sim_{num} K_W$ and such  that  the pull back of $G$ to $S'$ is $G'$. Let $t\colon W\to\bar{W}$ be the morphism onto the minimal model and let $E$ be the exceptional curve of $t$.  We have $GE=K_WE=-1$, hence $G=E+G_0$, where $G_0\ge 0$ and $G_0\sim_{num} t^*K_{\bar{W}}$. Assume that $W$ is properly elliptic and denote by $F$ a general fibre of the elliptic fibration of $W$. Then there is $\al\in \Q$, $\al>0$, such that $G_0\sim_{num}\al F$.
For $i=1,2$ let  $\Ga_i$ be the image of $Z_i$  in $W$. The curves $\Ga_1$ and $\Ga_2$ are $-4-$curves, hence we have $4=K_W(\Ga_1+\Ga_2)=E(\Ga_1+\Ga_2)+G_0(\Ga_1+\Ga_2)$. 
By construction, the curve $R_0$ does not meet the nodal curves $N_1,\dots, N_4$ of $W$ contained in the branch divisor $B_0$ of the double cover $V\to W$. Hence $EB_0=E(\Ga_1+\Ga_2)$ and  $G_0B_0=E(\Ga_1+\Ga_2)$ are both even. Moreover, we have $EB_0>0$, since otherwise $E$ would pull back on $S'$ to the disjoint union of two $-1-$curves, contradicting the minimality of $S'$, and $G_0B_0>0$, since otherwise $|F|$ would pull back on $S'$ to a pencil of elliptic curves, contradicting the fact that $S'$ is of general type. Hence we have $EB_0=G_0B_0=2$ and the pull back of $E$ on $S'$ is either a $-2-$curve or the union of two $-2-$curves meeting in a point. This is not possible, since $Z_1$ and $Z_2$ are the only $-2-$curves of $S'$ by construction. Hence we have reached a contradiction and the only possibility is that $W$ is birational to an Enriques surface.

\end{document}